\documentclass[11pt]{amsart}
\usepackage{amsfonts}
\usepackage{amssymb}
\usepackage{hyperref}
\usepackage{graphics}
\usepackage{hyperref}
\usepackage{amsmath}
\hypersetup{
colorlinks=true,
linkcolor=red,
citecolor=cyan,
}
\newtheorem{thm}{Theorem}[section]

\newtheorem{prop}[thm]{Proposition}

\newtheorem*{conj*}{Denjoy's Conjecture}

\theoremstyle{definition}

\newtheorem{defin}{Definition}

\newtheorem{exa}[thm]{Example}
\numberwithin{equation}{section}
\newcommand{\C}{\mathbb C}
\newcommand{\cN}{\mathcal N}
\newcommand{\cB}{\mathcal B}
\newcommand{\Scal}{\mathcal S}

\newcommand{\fall}[2]{{#1}^{\underline{#2}}}
\newcommand{\Lloc}{L^1_{{loc}}(\C)}

  \title[c. r g. and linear difference equation]{Completely Regular Growth of Transcendental Entire Solutions of Linear Difference Equations with Polynomial Coefficients }

    \author{Xuxu Xiang}
     \address{Xuxu Xiang\newline School of Mathematical Sciences, Guizhou Normal University, Guiyang, 550025, P.R. China. }
     \email{xiangxuxu@gznu.edu.cn}

\author{Jianren Long*}
     \address{Jianren Long \newline School of Mathematical Sciences, Guizhou Normal University, Guiyang, 550025, P.R. China. }
     \email{longjianren2004@163.com}

     \author{Mingjin Li}
     \address{Mingjin Li\newline School of Mathematical Sciences, Guizhou Normal University, Guiyang, 550025, P.R. China. }
     \email{limingjin2022@163.com}

   \keywords{Order of growth; entire function; completely regular growth; linear difference equation\\
   *Corresponding author.}

\begin{document}
\begin{abstract}Let \(f\) be a transcendental entire solution of order less than one of the linear difference equation
\[
\sum_{j=0}^{m}P_j(z)\Delta^j f(z)=0,
\qquad P_m\not\equiv0,
\]
where \(P_0,\ldots,P_m\) are polynomials. Applying the inverse Newton transform \(\mathcal S\), which maps the falling factorial \(z^{\underline n}\) to \(z^n\), we obtain a function \(g=\mathcal S f\) satisfying a nontrivial linear differential equation with polynomial coefficients. We further prove that the Newton transform preserves completely regular growth in the Levin–Pfluger sense. Consequently, \(f\) is of completely regular growth, and its zeros are regularly distributed. An example involving the  Gamma function shows that the conclusion may fail for solutions of order \(1\). Finally, we provide examples of linear difference equations with entire coefficients that admit solutions which are not of completely regular growth, even when all the coefficients are themselves of completely regular growth.

\end{abstract}

\maketitle

.

\section{Introduction and main results}\label{sec1}
 
For a meromorphic function $f$ on $\mathbb{C}$, let
$\Delta f(z)=f(z+1)-f(z)$ denote the   difference operator.
Set $\Delta^0f=f$ and define
$\Delta^n f=\Delta(\Delta^{n-1}f)$ for $n\geq 1$.
We consider the solutions of linear difference equation
\begin{equation}\label{eq1.1}
    \sum_{j=0}^{m}P_j(z)\Delta^j f(z)=0,
    \qquad P_m\not\equiv 0,
\end{equation}
where $P_0,\ldots,P_m$ are polynomials. If the coefficients \(P_0,\ldots,P_m\) are entire, then the linear difference equation \eqref{eq1.1} of order \(m\) possesses a fundamental system of \(m\) meromorphic solutions which are linearly independent over the field of \(1\)-periodic functions, see \cite{P1986}.

Using a difference analogue of Wiman--Valiron theory based on
binomial series, Ishizaki and Yanagihara
\cite[Theorem~1.1]{I2004} established the conclusion of
Theorem~1.1   below for transcendental entire solutions of order
$\rho(f)<1/2$. Chiang and Feng \cite[Theorem~7.3]{C2009}
extended the characterization of the order to the range
$\rho(f)<1$, showing that $\rho(f)$ is a positive rational
slope of the Newton polygon. They later refined their
difference Wiman--Valiron estimates and obtained the  following result.

\begin{thm}\cite[Theorem~4]{C2016}.
\label{tha}
Let \(f\) be an entire solution  of \eqref{eq1.1} with growth \(\rho(f) < 1\). Then  \(\rho(f)\) is a rational number which can be determined from a gradient of the corresponding Newton polygon of \eqref{eq1.1}. In particular,
\[
\log M(r, f) = L r^{\rho(f)} (1 + o(1)),
\]
where \(L > 0\), $\rho(f)=\limsup\limits_{r\to\infty}\frac{\log\log M(r,f)}{\log r}$,  \(M(r, f) = \max_{|z|=r} |f(z)|\).   
\end{thm}

Ishizaki and Wen \cite{I2025} established a criterion for
\eqref{eq1.1} to admit a transcendental entire solution of order less
than one by using binomial series. They \cite[Corollary 2.5.]{I2025} also showed that every transcendental entire solution
of \eqref{eq1.1} has positive order. 
Together with Theorem~\ref{tha}, this implies that
every transcendental entire solution \(f\) of
\eqref{eq1.1} with \(0<\rho(f)<1\) is of
\emph{perfectly regular growth} (p.r.g.).
Here, we say that an entire function \(f\) is of
p.r.g. if
$0<\rho(f)<\infty$
and its type satisfies
$ 
\tau(f):=\lim_{r\to\infty}
\frac{\log M(r,f)}{r^{\rho(f)}}\in(0,\infty).$
See \cite[Definition~7.4.1]{G2026} and
\cite[pp.~41--46]{V1949} for details.

An entire function $f$ of completely regular growth (c.r.g.)  implies perfectly p.r.g.see \cite[Theorem~7.5.7]{G2026}. Moreover, $f$ has c.r.g.  if and only if its zeros
are regularly distributed \cite[Theorem 1.2.1]{R1992}.
 The definitions of completely regular growth in the Levin--Pfluger sense and zeros
are regularly distributed  are stated as follows.

\begin{defin}[Completely regular growth {\cite[pp. 139-140]{J1980},\cite[Definition 7.5.1]{G2026}}] 
    We say that an entire function \(f\) of order \(\rho \in (0, \infty)\) and of mean type is of \textit{completely regular growth} (c.\,r.\,g.) if there exists a sequence of Euclidean discs \(D(z_k, r_k)\) satisfying
\begin{equation*}
    \sum_{|z_k| \le r} r_k = o(r) 
\end{equation*}
such that
\begin{equation*}
    \log|f(re^{i\theta})| = (h_f(\theta) + o(1))r^\rho, \quad re^{i\theta} \notin \bigcup_k D(z_k, r_k), 
\end{equation*}
as \(r \to \infty\) uniformly in \(\theta\). Here \(h_f(\theta) = \limsup_{r\to\infty} \frac{\log |f(re^{i\theta})|}{r^\rho}, \quad \theta \in \mathbb{R}\), is called the Phragm\'{e}n--Lindel\"{o}f indicator function.
\end{defin}

\begin{defin}[Angular
density and regularly distributed{\cite[pp. 6-7]{R1992}}]
    Let \(Z = \{w_k\}\) be an infinite sequence of complex numbers ordered by nondecreasing modulus such that \(|w_k| \to \infty\) as \(k \to \infty\). Suppose that the counting function \(n(r)\) of the points \(\{w_k\}\) has order \(\rho \in (0, \infty)\) and is of finite type with respect to \(\rho\). Let \(n_Z(r, \theta_1, \theta_2)\) denote the number of points \(\{w_k\}\) in the region
\[
\theta_1 < \arg z < \theta_2, \quad |z| \le r,
\]
each point being counted according to its multiplicity. We say that \(\{w_k\}\) has an \textit{angular density} with respect to \(\rho\) if for all but a countable set of values \(\theta_1, \theta_2\), the limit
\begin{equation*}
\lim_{r \to \infty} \frac{n_Z(r, \theta_1, \theta_2)}{r^\rho} 
\end{equation*}
exists as a finite number. If \(\rho \notin \mathbb{N}\), then the sequence \(\{w_k\}\) is called \textit{regularly distributed} if it has an angular density. For \(\rho \in \mathbb{N}\), the sequence \(\{w_k\}\) is called \textit{regularly distributed} if \(\{w_k\}\) has an angular density and, moreover, the limit
   $ \lim_{r\to\infty} \sum_{0<|w_k|\le r} w_k^{-\rho}$ 
exists as a finite number.
\end{defin}

Although Chiang and Feng  referred to the growth property
obtained in Theorem  \ref{tha} as completely regular growth, their
result actually establishes perfectly regular growth (p.r.g.)
in the sense adopted here.
This is different from completely regular growth (c.r.g.)
in the Levin--Pfluger sense considered in this paper, which
requires the corresponding angular asymptotic behavior
outside an exceptional set of discs.
According to the discussion above,
a natural question: does every transcendental entire solution $f$
of \eqref{eq1.1} with $\rho(f)<1$ have c.r.g.? In this paper, we resolve this question by establishing the following result.

\begin{thm}
\label{th1.2}
    Let \(f\) be an  transcendental entire solution  of \eqref{eq1.1} with  \(\rho(f) < 1\). Then $f$ is of c.r.g. and  the zeros of \(f\) are regularly distributed.
\end{thm}
 
The following example involving the \(\Gamma\) function
shows that the restriction \(\rho(f)<1\) is essential
in general. Some examples demonstrating that the condition of polynomial coefficients cannot be relaxed to transcendental entire coefficients will be given in Section~4.
\begin{exa}
\label{e1.3}
Let \(f(z) = \frac{1}{\Gamma(z)}\), which is an entire function. Since the Gamma function satisfies \(z\Gamma(z) = \Gamma(z+1)\), we deduce that \(f\) satisfies the difference equation \(z\Delta f(z) + (z - 1)f(z) = 0\). However, it is known that \(\rho(f) = 1\) and \(\tau(f) = \infty\) (see \cite[pp. 307 and pp.309]{M1965}), which implies that \(f\) is not of c.\,r.\,g.
\end{exa}

In what follows, we briefly summarize the main ideas of the proof and the organization of this paper. The starting point of our argument is the binomial-series approach of Ishizaki and Wen \cite{I2021}. They \cite[Corollary A.1.]{I2021} gave the   following representation for
transcendental entire solution  of \eqref{eq1.1} with  \(\rho(f) < 1\):

\begin{equation}\label{eq1.2}
 f(z)=\sum_{n=0}^{\infty}a_n\fall{z}{n},
 \qquad a_n=\frac{\Delta^nf(0)}{n!}.
\end{equation}
Here
   \(z^{\underline{0}} = 1\) and
\begin{equation}
    z^{\underline{n}} = z(z-1)\cdots(z-n+1) = n!\binom{z}{n}, \quad n=1,2,3,\dots,
\end{equation}
which is called a \textit{falling factorial}.

In Section~2, we will give the definitions of the Newton transform $\cN$
 and its inverse $\Scal$.  
$\Scal$
 is the 
map that sends $\fall{z}{n}$ to $z^n$.  We call $\Scal$ the inverse Newton transform, in reference to the transform
$\cN(z^n)=\fall{z}{n}$, called the Newton transform in
\cite[Example~5.3, Eq.~(5.3)]{c2023}.
Since $\Delta\fall{z}{n}=n\fall{z}{n-1}$, we have
\[
    \Scal(\Delta f)(z)
    =\Scal\!\left(\sum_{n=1}^{\infty}na_n\fall{z}{n-1}\right)
    =\sum_{n=1}^{\infty}na_nz^{n-1}
    =Dg(z),
\]
   where  $g(z)=\sum_{n=0}^{\infty}a_nz^n$ is the ordinary power
series  and     $D=\frac{d}{dz}$.  We will show that $g=\Scal(f)$ is an entire function of the same order
as $f$ and establish how $\Scal$ acts on  $P_j(z)\Delta^jf(z)$.  These identities allow us to transform
\eqref{eq1.1} into a linear differential equation with polynomial
coefficients for $g$.
  It is well known that transcendental entire solutions of such equations possess c.\,r.\,g. \cite[pp. 104--112]{P1984}\cite{S1989}. Since $f=\cN(g)$, we are led to prove that the Newton transform $\cN$ preserves the c.\,r.\,g. of $g$.
  
  Indeed, inspired by the works of \cite{Beise2011,Levin1996}, we formulate the Newton transform in the form of the Borel transform. 
The potential-theoretic form of c.\,r.\,g. for \(g\) is equivalent to the \(L_{\text{loc}}^1\)-convergence of \(\frac{\log|g(Rz)|}{R^\rho}\) to \(H_g(z) = |z|^\rho h_g(\theta)\) as \(R \to \infty\) \cite{Azarin2009,Eremenko1998,Levin1996}. This means that for every compact set \(K \Subset \C\), 
\begin{equation}\label{eq:CRG}
\int_K \left|\frac {\log |g(Rz)|}{R^{\rho}}-H_g(z)\right|dA\to 0 \quad \text{as } R \to \infty.
\end{equation}
Combining the Borel integral representation with
Hartogs' compactness principle for subharmonic
functions, we prove that the normalized logarithmic
potentials of \(g\) and its Newton transform
\(\cN(g)\) have the same limit in
\(L_{\mathrm{loc}}^1(\mathbb C)\).
Consequently, the Newton transform preserves
completely regular growth and the indicator function.
We establish this preservation property in Section 3.

In Section 4, we complete the proof of
Theorem~\ref{th1.2} by applying the results on
binomial series and the inverse Newton transform
established in Section 2, together with the
preservation property proved in Section 3.
We then provide examples to show that the conclusion
of Theorem~\ref{th1.2} may fail when the polynomial
coefficients are replaced by transcendental entire
coefficients, even if these coefficients are of c.r.g.

\section{Binomial Series and the Inverse Newton Transform}

In this section, we   present some results on the binomial series and the Inverse Newton transform. Let $\{b_n\}$ be a   sequence.
Define
\begin{align}
\label{2.1}
\chi(\{b_n\})=\limsup_{n\to\infty}
          \frac{n\log n}{-\log|b_n|}, 
\end{align}
where the expression $\limsup_{n\to\infty}
          \frac{n\log n}{-\log|b_n|}$ is considered to be zero whenever $b_n=0$.

Consider the binomial series of the form  
\begin{align}
\label{eq2.2}
   Y(z)=\sum_{n=0}^{\infty}b_n\fall{z}{n}.   
\end{align}
  Ishizaki and Wen \cite{I2021} used \eqref{2.1} to study the growth of $Y(z)$.
\begin{thm}\cite[Theorem 1.1]{I2021}
\label{th2.1}
Suppose that \(\chi(\{b_n\}) < 1\). Then, the formal series \(Y(z)\) given by \eqref{eq2.2} converges to \(y(z)\) uniformly on every compact subset in \(\mathbb{C}\). Moreover, the order of growth for \(y(z)\) coincides with \(\chi(\{b_n\})\).
\end{thm}
Their \cite[Corollary A.1.]{I2021} gives an   entire $Y$  of  $\rho(Y)<1$  can be represented by a binomial series. Together with Theorem \ref{th2.1}, we obtain the following result.

\begin{thm}
\label{th2.2}
\(Y(z)\) given by \eqref{eq2.2} with
$\chi(\{b_n\}) < 1$ $\iff $  $Y$ is an entire function of  $\rho(Y)<1$.  Consequently, $\rho(Y)=\chi(\{b_n\})$.
\end{thm}
\begin{proof}[proof of Theorem \ref{th2.2}]

\textbf{Sufficiency.} Since \(\chi(\{b\}) < 1\), Theorem \ref{th2.1} guarantees that \(Y\) converges uniformly on every compact subset of \(\mathbb{C}\) to an entire function, which we also denote by \(Y\). Thus, \(\rho(Y) = \chi(\{b\}) < 1\).
 
 \textbf{Necessity}. Since $\rho(Y)<1$, then \cite[Corollary A.1.]{I2021} implies $Y$ can be represented
by a binomial series \eqref{eq2.2} and $  b_n=\frac{\Delta^nY(0)}{n!}.$  Next, we estimate \(b_n\) by express the difference as an integral of the derivative.

 Repeated this integration \[
\Delta Y(z) = Y(z + 1) - Y(z) = \int_0^1 Y'(z + t) \, dt 
\]gives
\[
\Delta^n Y(0) = \int_{[0,1]^n} Y^{(n)}(t_1 + \cdots + t_n) \, dt_1 \cdots dt_n.
\]
For any \(x \in [0, n]\), applying Cauchy's derivative estimate on a circle centered at \(x\) with radius \(R\), we get
\[
|Y^{(n)}(x)| \le \frac{n!}{R^n} M(n + R, Y).
\]
This gives
\[
|b_n| = \frac{|\Delta^n Y(0)|}{n!} \le \frac{M(n + R, Y)}{R^n}.
\]

Take \(\rho(Y) < \sigma < 1\). By the definition of the order of growth, we obtain
\[
|Y(z)| \le C_\sigma e^{|z|^\sigma},
\]
where $C_\sigma$ is a nonzero constant.  
Choosing \(R = n^{1/\sigma}\), since \(\sigma < 1\), we have \(R \ge n\), and therefore
$(n + R)^\sigma \le (2R)^\sigma = 2^\sigma n.$ Thus
\begin{align}
\label{eq2.3}
   |b_n| \le C_\sigma e^{2^\sigma n} n^{-n/\sigma}. 
\end{align}
Therefore,
\[
-\log|b_n| \ge \frac{n}{\sigma}\log n - (2^\sigma n + \log C_\sigma), \qquad \chi(\{b\}) = \limsup_{n\to\infty} \frac{n\log n}{-\log|b_n|} \le \sigma < 1.
\]
Since \(\chi(\{b\}) < 1\), Theorem \ref{th2.1} guarantees that $\rho(Y)=\chi(\{b_n\})$.
\end{proof}
\begin{defin}[Newton transform and inverse Newton transform]
\label{def:newton-transforms}
For $F$ with Taylor expansion
$F(z)=\sum_{n=0}^{\infty}a_n z^n$, define its
\emph{Newton transform} by
\[
    (\cN F)(z)
    :=\sum_{n=0}^{\infty}a_n\fall{z}{n}.
\]
For entire function $Y$ given by \eqref{eq2.2}, define its \emph{inverse Newton transform} by
\[
    (\Scal Y)(z)
    :=\sum_{n=0}^{\infty}b_nz^n.
\]
\end{defin}
 We give some   propositions of   inverse Newton transform $\Scal$. 
 \begin{prop}
 \label{p1}
     If \(Y\) is a transcendental entire function of $\rho(Y)<1$, then \[ F :=   \Scal\left( Y=\sum_{n=0}^{\infty} b_n \fall{z}{n} \right) = \sum_{n=0}^{\infty} b_n z^n\]
 is an entire function and $\rho(Y)=\rho(F)$,  $b_n=\frac{\Delta^nY(0)}{n!} $. Moreover, \(\Scal\) is bijective and \(\Scal^{-1}=\cN\)
on the class of transcendental entire functions of order less than one.
 \end{prop}

 \begin{proof}[Proof of Proposition \ref{p1}]
  Since $\rho(Y)<1$, then \cite[Corollary A.1.]{I2021} implies $Y$ can be represented
by a binomial series \eqref{eq2.2} and $  b_n=\frac{\Delta^nY(0)}{n!}.$ 
By an argument analogous to the proof of Theorem \ref{th2.2}, we can obtain estimate \eqref{eq2.3} of $b_n$,
\[
|b_n|^{1/n} \le C_\sigma^{1/n} e^{2^\sigma} n^{-1/\sigma} \longrightarrow 0.
\]
Therefore,  \(F\) is an entire function. By  Taylor coefficient formula for the order of an
entire function \cite[Theorem 2.2.2]{B1954},    then we obtain  $\rho(F)=\chi(\{b\})$. Together this with  Theorem \ref{th2.2} gives 
$\rho(Y)=\chi(\{b_n\})=\rho(F)$.

 Since
   $ \cN(\Scal(Y))
    = \sum_{n=0}^{\infty} b_n \fall{z}{n}
    = Y.$ 
Thus \(\Scal\) is injective.
Next we prove surjectivity, let \(F_1(z)=\sum_{n=0}^{\infty} a_n z^n\)
be a transcendental entire function of order less than one.
By the Taylor coefficient formula for the order of an entire function \cite[Theorem 2.2.2]{B1954},
\(\chi(\{a_n\})=\rho(F_1)<1\). It follows from
Theorem~\ref{th2.2} that
\[
    Y_1(z):=\cN(F_1)(z)
    =\sum_{n=0}^{\infty} a_n \fall{z}{n}
\]
is an entire function of  $\rho(Y_1)=\chi(\{a_n\})<1$.
Thus,
\[
    \Scal(\cN(F_1))(z)=\Scal(Y_1)(z)
    =\sum_{n=0}^{\infty} a_n z^n=F_1(z).
\]
Therefore \(\Scal\) is bijective and \(\Scal^{-1}=\cN\)
on the class of transcendental entire functions of order less than one.
\end{proof}

\begin{prop}
\label{p2.4}
    Let $Y$ be an transcendental entire function with $\rho(Y)<1$ and  let $P(z)$ be a polynomial. Then, $\Scal(Y)=F$, $\Scal (P(z)\Delta^j Y)=P(A)D^jF$, where $A=z(I+D)$ is differential operator with polynomial coefficients, $F=\sum_{n=0}^{\infty} b_n z^n$
 is an entire function and    $b_n=\frac{\Delta^nY(0)}{n!}$. 
\end{prop}

\begin{proof}[proof of Proposition \ref{p2.4}]
 Since $\rho(Y)<1$, then \cite[Corollary A.1.]{I2021} implies $Y$ can be represented
by a binomial series \eqref{eq2.2} and $  b_n=\frac{\Delta^nY(0)}{n!}.$  Using Proposition \ref{p1} implies $\Scal(Y)=F=\sum_{n=0}^{\infty} b_n z^n$
 is an entire function.  By $\Delta\fall{z}{n}=n\fall{z}{n-1}$, we have  
\[
    \Scal(\Delta Y)(z)
    =\Scal\!\left(\sum_{n=1}^{\infty}nb_n\fall{z}{n-1}\right)
    =\sum_{n=1}^{\infty}nb_nz^{n-1}
    =DF(z).
\]

Ishizaki and Wen \cite[Lemma 4.1]{I2021} gave  that \[
zY(z) = \sum_{n=1}^{\infty} (n b_n + b_{n-1})\fall{z}{n}, \qquad z\Delta Y(z) = \sum_{n=1}^{\infty} n \big( (n+1)b_{n+1} + b_n \big) \fall{z}{n}
\]. Therefore, \[\Scal(zY)=\sum_{n=1}^{\infty} (n b_n + b_{n-1})z^n=zD\Scal(Y)+z\Scal(Y)=A\Scal(Y).\]
Since $z^kY=z(z^{k-1}Y)$, then $\Scal(z^kY)=A^k\Scal(Y)$, $\Scal(P(z)Y)=P(A)\Scal(Y)$  and $\Scal(P(z)\Delta^jY)=P(A)\Scal(\Delta^jY)=P(A)D^jF$.

\end{proof}

\section{Preservation of Completely Regular Growth under the Newton Transform}

An entire function \(f\) is said to be of exponential type if \(\rho(f) < 1\), or if \(\rho(f) = 1\) and \(\tau(f) < \infty\)\cite[pp.7]{B1954}.
  The  function \(f(z) = \sum_n a_n z^n\) is an entire function  of exponential
type if and only if
\begin{equation}\label{eq:Borel}
\mathcal{B}f(z) := \sum_{n=0}^{\infty} f^{(n)}(0)z^{-n-1} = \sum_{n=0}^{\infty} n! a_n z^{-n-1} 
\end{equation}
 is convergent for some (finite) $z$ \cite[5.3. The indicator diagram]{B1954}.
$\mathcal{B}f(z)$ is called the Borel transform of \(f\). P\'olya's theorem (See \cite[Theorem 2, pp.66]{Levin1996} for details) yields 
\begin{equation}\label{eq:Polya}
f(z) = \frac{1}{2\pi i} \int_{|\xi|=r} \mathcal{B}f(\xi) e^{\xi z} d\xi. 
\end{equation}
Now we give the main result in this section.
\begin{thm}
\label{thm3.1}
Let $f$ and $g$ be entire functions of the same order $\rho\in(0,1)$
and of mean type, with $f=\cN(g)$ and $g=\Scal(f)$.
If $g$ is of c.\,r.\,g., then $f$ is of c.\,r.\,g. at the same order and $h_f=h_g$. 
\end{thm}

\begin{proof}

We use some ideas from \cite{Azarin2009,Levin1996,Beise2011}. Choose $C,B>0$ such that
\begin{equation}
\label{eq:newton-crg-type}
  |u(z)|\le C e^{B|z|^\rho},\qquad z\in\C,\quad u\in\{f,g\}.
\end{equation} 
Using the Cauchy's derivative estimate on
the circle of radius $k^{1/\rho}$ centered at $w$. By
\eqref{eq:newton-crg-type}, the fact \(0<\rho <1\) and $k!\le k^k$, for $k\ge1$ we obtain
\begin{align}\label{eq:newton-crg-derivative}
\begin{aligned}
    |u^{(k)}(w)|
  &\le C\frac{k!}{k^{k/\rho}}
        \exp\!\bigl(B(|w|+k^{1/\rho})^\rho\bigr)\\
  &\le C e^{B|w|^\rho}(e^B k^{-\alpha})^k,
     \qquad \alpha=\frac{1-\rho}{\rho}>0. 
\end{aligned}
 \end{align}

For \(g(w)=\sum_{n=0}^\infty\frac{g^{(n)}(0)}{n!}w^n\) is exponential type, choosing \(|s|<1\) such that \((1+\zeta)^w\) is well defined and the summation by terms yields $$
\begin{aligned}
    \cN(g)(w)=&\sum_{n=0}^\infty\frac{g^{(n)}(0)}{n!}\fall{z}{n}=\sum_{n=0}^\infty {g^{(n)}(0)}\binom{w}{n}\\
    &=\frac{1}{2\pi i}\int_{|\zeta|=s} (\cB g)(\zeta)e^{w\log(1+\zeta)}d\zeta.
\end{aligned}$$
 By induction, \(\Delta^jf(a)=\frac{1}{2\pi i}\int_{|\zeta|=s}(\cB f)(\zeta)e^{a\zeta}(e^\zeta-1)^jd\zeta\), hence \(\Delta^j f(0)=\frac{1}{2\pi i}\int_{|\zeta|=s}(\cB f)(\zeta)(e^\zeta-1)^jd\zeta\).\\
For \(f(w)=\sum_{j=0}^\infty\frac{\Delta^jf(0)}{j!}\fall{w}{j}\),  summation by terms yields
$$
\begin{aligned}
    \Scal f(w)=\sum_{j=0}^\infty\frac{\Delta^jf(0)}{j!}w^j=\frac{1}{2\pi i}\int_{|\zeta|=s} \cB f(\zeta)e^{w(e^{\zeta}-1)}d\zeta.
\end{aligned}
$$
Define \[
    (T,u)=:(\cN,g)\,\,\text{or}\,\, (\Scal, f).
\]
For \((T,u)=:(\cN,g)\,\,\text{or}\,\, (\Scal, f)\), set
\[
    \phi_T(\zeta)=\begin{cases}
     \log(1+\zeta),\,\, &T=\cN,\\
   e^{\zeta}-1,\,\, &T=\Scal.
    \end{cases}
\]
The above argument yields 
\[
    Tu(w)=\frac{1}{2\pi i}\int_{|\zeta|=s} \cB u(\zeta)e^{w\phi_T(\zeta)}d\zeta.
\]
Noticing that the Borel transform of the function \(\zeta\mapsto u(w+\zeta)\) is \(\cB_w(\zeta)=\sum_{n=0}^\infty\frac{u^{(n)}(w)}{\zeta^{n+1}}\), differentiating term by term shows that 
\[
    \partial_w B_w u(\zeta)=\zeta B_wu(\zeta)-u(w),
\]
hence integrating term by term shows that
\[
\cB_wu(\zeta)-e^{w\zeta}\cB u(\zeta)=-\int_0^we^{(w-z)\zeta}u(z)dz.
\]
Set \(b_T(\zeta)=\phi_T(\zeta)-\zeta\), then \(b_T(\zeta)\) is holomorphic and  \(|b_T(\zeta)|=O(|\zeta|^2)\). Cauchy's integration theorem yields 
\[
    \frac{1}{2\pi i}\int_{|\zeta|=s}\cB_wu(\zeta)e^{wb_T(\zeta)}-e^{w\phi_T(\zeta)}\cB u(\zeta)d\zeta=0,
\]
thus, 
\[
   Tu(w)= \frac{1}{2\pi i}\int_{|\zeta|=s}\cB_wu(\zeta)e^{wb_T(\zeta)}d\zeta.
\]

By \cite[equation 4, pp.66]{Levin1996},  for \(\zeta=se^{i\theta}\),
\[
    B_w(se^{i\theta})=\int_0^\infty u(w+xe^{-i\theta})e^{-i\theta}e^{-sx}dx,
\]
hence 
\( |Tu(w)|\leq I_R+E_R\), where 
\[
    I_R=\frac{1}{2\pi}\int_0^{2\pi}\left|e^{wb_T(se^{i\theta})}\right|s\int_0^{\delta R}|u(w+xe^{-i\theta})|e^{-sx}dxd\theta,
\]  
\[
    E_R=\frac{1}{2\pi}\int_0^{2\pi}\left|e^{wb_T(se^{i\theta})}\right|s\int_{\delta R}^{\infty}|u(w+xe^{-i\theta})|e^{-sx}dxd\theta.
\]
Fixed \(L, \delta>0\), set \(M_u(w,\delta R)=\sup_{|\zeta-w|\leq \delta R}|u(\zeta)|\), choosing \(s=s_R=R^{\rho-1}\log R\),  suppose \(|w|\leq \delta R\), then for \(x\leq \delta R\), 
\[
    |I_R|\leq CM_u(w,\delta R)e^{CLR^{2\rho-1}(\log R)^2}C\leq e^{o(R^\rho)}M_u(w,\delta R).
\]
For \(x\geq\delta R\), by \eqref{eq:newton-crg-type} and \(0<\rho<1\), then \(|u(w+xe^{-i\vartheta})|
\le C e^{AL^\rho R^\rho+Ax^\rho}\) and notice that \(
\frac{Ax^\rho}{s_Rx}
\le \frac{A\delta^{\rho-1}}{\log R}.
\to 0\), hence 
\[s_R\int_{\delta R}^\infty |u(w+xe^{-i\theta})e^{-s_Rx}dx\leq Ce^{AL^\rho R^\rho}s\int_{\delta R}^\infty e^{Ax^\rho-s_Rx}dx\leq
2Ce^{AL^\rho R^\rho-\frac{\delta}{2}s_RR}.
\] 
Noticing again 
\[
    \frac{\log2C+C_LRs_R^2+AL^\rho R^\rho}{Rs_R}\to 0, 
\]
it follows that \[
\begin{aligned}
   E_R&\leq 
2Ce^{C_LRs_R^2+AL^\rho R^\rho-\frac{\delta}{2}s_RR}\\
&=e^{\log2C+C_LRs_R^2+AL^\rho R^\rho-\frac{\delta}{2}s_RR}\leq e^{-\frac{\delta}{4}RS_R}=e^{-\frac{\delta}{4}R^\rho\log R}.  
\end{aligned}
   \]
Thus 
\begin{equation}\label{eq:module}
    |Tu(w)|\leq e^{o(R^\rho)}M_u(w,\delta R)+e^{-cR^\rho\log R}.
\end{equation}

Set \[
U_R(z) = \frac {1}{R^\rho} \log |g(Rz)|, \qquad V_R(z) = \frac {1}{R^\rho} \log |f(Rz)|.
\]
For \(|z|\leq L\), by \eqref{eq:newton-crg-type}, we have
\[
V_{R_j}(z)
=\frac{\log|f(R_jz)|}{R_j^\rho}
\leq \frac{\log C}{R_j^\rho}+B|z|^\rho
\leq C_L,
\]
where \(C_L\) is a positive constant depending on \(L\).

Choosing \(\xi\) such that \(f(\xi)\neq 0\), set
\(z_j=\frac{\xi}{R_j}\). Then
\[
V_{R_j}(z_j)
=\frac{\log|f(\xi)|}{R_j^\rho}\to 0.
\]
Since \(z_j\to 0\), these points eventually lie in the closed unit
disk. Hence, \(V_{R_j}\) cannot converge uniformly to
\(-\infty\) on the closed unit disk. Since \(V_{R_j}\) is
subharmonic and locally uniformly bounded above, by
\cite[Theorem~4.1.9]{Hormander1990}, after passing to a
subsequence, we have
\[
\|V_{R_j}-V\|_{L^1_{\mathrm{loc}}(\C)}\to 0,
\]
where \(V\) is a subharmonic function. By \eqref{eq:CRG},
\(U_{R_j}\) also converges to \(H_g\) in
\(L^1_{\mathrm{loc}}(\C)\). By Riesz's theorem, there exists
a common subsequence, still denoted by \(\{R_j\}\), such that
\[
U_{R_j}(z)\to H_g(z),\qquad
V_{R_j}(z)\to V(z)
\quad\text{a.e.},\quad j\to\infty,
\]
where \(H_g(re^{i\theta})=r^\rho h_g(\theta)\).

Let \(f=\cN g\), \(w=Rz\), and \(\zeta=R\eta\). Then
\[
\frac{\log M_g(Rz,\delta R)}{R^\rho}
=\sup_{|\eta-z|\leq\delta}U_R(\eta).
\]
By \eqref{eq:module}, we obtain
\begin{equation}\label{eq:VR}
V_R(z)\leq \frac{\log 2}{R^\rho}
+\max\left\{
\sup_{|\eta-z|\leq\delta}U_R(\eta)+o(1),
-c\log R
\right\}.
\end{equation}
Similarly,
\begin{equation}\label{eq:UR}
U_R(z)\leq \frac{\log 2}{R^\rho}
+\max\left\{
\sup_{|\eta-z|\leq\delta}V_R(\eta)+o(1),
-c'\log R
\right\}.
\end{equation}

For fixed \(\delta>0\) and \(z\in\C\), by the subharmonicity
of \(V_{R_j}\), for \(|\zeta-z|\leq\delta\), we have
\[
\begin{aligned}
V_{R_j}(\zeta)
&\leq \frac{1}{\pi\delta^2}
\int_{D(\zeta,\delta)}V_{R_j}(\omega)\,dA(\omega)\\
&\leq \frac{1}{\pi\delta^2}
\int_{D(\zeta,\delta)}V(\omega)\,dA(\omega)
+\frac{\|V_{R_j}-V\|_{L^1(D(z,2\delta))}}
{\pi\delta^2}\\
&\leq \sup_{|\eta-z|\leq 2\delta}V(\eta)
+\frac{\|V_{R_j}-V\|_{L^1(D(z,2\delta))}}
{\pi\delta^2}.
\end{aligned}
\]
Thus,
\[
\limsup_{j\to\infty}
\sup_{|\zeta-z|\leq\delta}V_{R_j}(\zeta)
\leq \sup_{|\eta-z|\leq 2\delta}V(\eta).
\]
Similarly, since \(U_{R_j}\) is subharmonic, we obtain
\[
\limsup_{j\to\infty}
\sup_{|\zeta-z|\leq\delta}U_{R_j}(\zeta)
\leq \sup_{|\eta-z|\leq 2\delta}H_g(\eta).
\]

By \eqref{eq:VR}, noting that
\(-c\log R_j\to-\infty\), we have
\[
\begin{aligned}
V(z)
&=\lim_{j\to\infty}V_{R_j}(z)\\
&\leq \limsup_{j\to\infty}
\sup_{|\zeta-z|\leq\delta}U_{R_j}(\zeta)\\
&\leq \sup_{|\eta-z|\leq 2\delta}H_g(\eta)
\end{aligned}
\]
for almost every \(z\in\C\). Similarly, by \eqref{eq:UR},
\[
\begin{aligned}
H_g(z)
&=\lim_{j\to\infty}U_{R_j}(z)\\
&\leq \limsup_{j\to\infty}
\sup_{|\zeta-z|\leq\delta}V_{R_j}(\zeta)\\
&\leq \sup_{|\eta-z|\leq 2\delta}V(\eta)
\end{aligned}
\]
for almost every \(z\in\C\).

Taking \(\delta=1/n\), \(n=1,2,\ldots\), and letting
\(n\to\infty\), since \(H_g\) is continuous and \(V\) is
upper semicontinuous, we obtain
\[
V(z)\leq H_g(z),\qquad H_g(z)\leq V(z)
\]
for almost every \(z\in\C\). Hence,
\[
V(z)=H_g(z)\quad\text{a.e. in }\C.
\]
Consequently,
\begin{equation}\label{eq:VRJ}
\|V_{R_j}-H_g\|_{\Lloc}\to 0,
\qquad j\to\infty.
\end{equation}

Suppose  that \(\{V_R\}\) does not converge
to \(H_g\) in \(L^1_{\mathrm{loc}}(\C)\). Then there exist
\(\varepsilon_0>0\), \(A>0\) and a sequence \(R_j\to\infty\)
such that
\[
\int_{D(0,A)}|V_{R_j}-H_g|\,dA\geq\varepsilon_0
\]
for all \(j\). Applying the above argument to this sequence,
we can choose a subsequence \(\{R_{j_k}\}\) such that
\[
\|V_{R_{j_k}}-H_g\|_{L^1(D(0,A))}\to 0,
\qquad k\to\infty.
\]
This is a contradiction. Therefore,
\[
\|V_R-H_g\|_{\Lloc}\to 0,
\qquad R\to\infty.
\]
This completes the proof.
\end{proof}

\section{ Proof of Theorem \ref{th1.2} and the Example of Transcendental Entire Coefficients}
\begin{proof}[Proof of the   Theorem \ref{th1.2}]

Let $f$ be an transcendental entire solution of  \eqref{eq1.1}. Then \cite[Corollary 2.5.]{I2025} gives $\rho(f)>0$. 
By \cite[Corollary A.1]{I2021}, we have
\[
    f(z)=\sum_{n=0}^{\infty}a_n\fall{z}{n},
    \qquad a_n=\frac{\Delta^n f(0)}{n!}.
\]
By Proposition~\ref{p1},
\[
    g(z):=\Scal(f)(z)=\sum_{n=0}^{\infty}a_nz^n
\]
is an entire function of order \(\rho(g)=\rho(f)<1\), and
\(\cN(g)=f\).

Let \(k\)
be the least index for which \(P_k\not\equiv0\). Since \(\Scal\) preserves the degree and
leading coefficient of every nonzero polynomial, we have
\(\Scal(P_k)\not\equiv0\), 
\[
    \sum_{j=0}^{m}P_j(A)D^j(z^k)
    =\Scal\!\left(\sum_{j=0}^{m}P_j(z)\Delta^j\fall{z}{k}\right)
=\Scal(k!P_k)=k!\,\Scal(P_k)\not\equiv0,
\]
Thus $\sum_{j=0}^{m}P_j(A)D^j\not=0$.
Since \(f\) solves \eqref{eq1.1} and \(\Scal(0)=0\) and by linearity and Proposition~\ref{p2.4}, it follows that
\[
    \Scal\!\left(\sum_{j=0}^{m}P_j(z)\Delta^j f(z)\right)
    = \sum_{j=0}^{m}P_j(A)D^j g(z)=0.
\]
Hence \(g\) is a transcendental entire solution of a nontrivial
homogeneous linear differential equation with polynomial coefficients.
By the classical results in \cite[pp.~104--112]{P1984} and \cite{S1989},
\(g\)  is of c.r.g. Thus,  Theorem \ref{thm3.1} impels  $f$ is of c.r.g.  This shows that the zeros of $f$ are regularly distributed \cite[Theorem 1.2.1]{R1992}.
\end{proof}

  Gol'dberg and Ostrovski\u{\i} \cite[p.~300]{H1994} asked \emph{whether a transcendental entire solution \(f\) of \(f^{(n)} + a_{n-1}f^{(n-1)} + \dots + a_1f' + a_0f = 0\) of finite order must be of c.\,r.\,g. if the coefficients \(a_0, \dots, a_{n-1}\) are of c.\,r.\,g.} Without the condition of c.\,r.\,g. for the coefficients, the answer to this
question is negative, see \cite[p.~300]{H1994}. For some further discussion and results concerning this question, we refer to \cite{B2020,H2015,Wen2018}. In particular, Bergweiler \cite{B2023} answered this question in the negative.

 A natural question is to consider the c.\,r.\,g. of transcendental entire solutions of linear difference equations \eqref{eq1.1} when the coefficients are entire functions, as well as the difference analogue of the Gol'dberg and Ostrovski\u{\i} question: \emph{whether a transcendental entire solution \(f\) of \eqref{eq1.1} of finite order must be of c.\,r.\,g. if the coefficients are of c. r. g}. However, we provide the following example to show that there exists a solution of a difference equation with entire coefficients which is not of c.\,r.\,g., even if the coefficients are of  c.\,r.\,g.

\begin{exa} Take
\[
F(z) = \prod_{k=1}^{\infty} \left(1 - \frac{z}{4^k}\right)^{2^k}.
\]
Since \(\sum_k 2^k/4^k < \infty\), the product converges locally uniformly. Let \(n_F(r)\) be its zero-counting function, counting multiplicities. Then \(n_F(t) \le 2\sqrt{t}\), so
\[
\log M(r,F) = \sum_{k \ge 1} 2^k \log \left(1 + \frac{r}{4^k}\right) = \int_0^{\infty} \frac{r n_F(t)}{t(t+r)} \, dt \le 2\pi\sqrt{r}.
\]
On the other hand, \(\log M(4^k, F) \ge 2^k \log 2\). Therefore, \(\rho(F) = 1/2\), and its type is finite and non-zero.
However,
\[
\frac{n_F(4^k)}{\sqrt{4^k}} = \frac{2^{k+1} - 2}{2^k} \longrightarrow 2, \qquad \frac{n_F(2 \cdot 4^k)}{\sqrt{2 \cdot 4^k}} \longrightarrow \sqrt{2}.
\]
By \cite[Chapter~III, \S3, Theorem~3]{J1980},
c.\,r.\,g. with respect to $r^{1/2}$
would imply the existence of
$\lim_{r\to\infty}n_F(r)/r^{1/2}$.
Since this limit does not exist, $F$ is not of c.\,r.\,g. with respect to $r^{1/2}$.

(1) Consider
\[
F(z)\Delta f(z) + (F(z) - F(z+1))f(z) = 0.
\]
Both coefficients are entire functions. Direct substitution shows that \(f = F\) is a solution, since
\[
F(z)(F(z+1) - F(z)) + (F(z) - F(z+1))F(z) = 0.
\]
But $F$ does not have c.\,r.\,g.

(2) Consider
\[
 e^z F(z)\Delta f(z) + e^z (F(z) - F(z+1))f(z) = 0. 
\]
It still has \(f = F\) as a solution. However, the coefficients are of c.~r.~g. Since for any non-zero entire function \(H\), if \(\rho(H) < 1\), then \(e^z H(z)\) is an entire function of c.~r.~g. of order \(1\) with indicator \(\cos \theta\).  
\end{exa}

\begin{exa}
Consider
\[
 ze^z \Delta f(z) + (z - 1)e^z f(z) = 0. 
\]
Both coefficients are entire functions of c.\,r.\,g., 
while the   solution \(f = 1/\Gamma\)   does not have c.\,r.\,g.
\end{exa}

	\section*{Declarations}
	\begin{itemize}
		
		\item\noindent{\bf Funding}
 Xuxu Xiang was supported by Graduate Research Fund Project of Guizhou Province (2025YJSKYJJ107). Jianren Long was supported by the National Natural Science Foundation of China (Grant No. 12261023, 11861023) 
		\item \noindent{\bf Conflicts of Interest}
		The authors declare that there are no conflicts of interest regarding the publication of this paper.

        \item\noindent{\bf Author Contributions}
All authors contributed to the study conception and design. All authors read and approved the final manuscript. 

 \item\noindent{\bf Acknowledgments} The first author wishes to express his sincere gratitude to Professor Zhitao Wen of Shantou University for his kind help and valuable discussions on linear difference equations during his visit to Shantou University in 2025.
 
	\end{itemize}


\begin{thebibliography}{99}
	\normalsize
	\baselineskip=17pt

 \bibitem{Azarin2009}
V.~S.~Azarin,
\emph{Growth Theory of Subharmonic Functions},
Birkh\"auser, Basel, 2009.
\url{https://doi.org/10.1007/978-3-7643-8886-7}.

\bibitem{B2020}
A.~I.~Bandura and O.~B.~Skaskiv,
Boundedness of the \(l\)-index and the completely regular growth
of entire functions (Ukrainian),
\emph{Ukrain. Mat. Zh.} \textbf{72} (2020), no.~3, 316--325;
English translation:
\emph{Ukrainian Math. J.} \textbf{72} (2020), no.~3, 358--369.
\url{https://doi.org/10.1007/s11253-020-01787-w}.

\bibitem{Beise2011}
H.-P.~Beise and J.~M\"uller,
Growth of (frequently) hypercyclic functions for differential operators,
\emph{Studia Math.} \textbf{207} (2011), no.~2, 97--115.
\url{https://doi.org/10.4064/sm207-2-1}.

\bibitem{B2023}
W.~Bergweiler,
A question of Gol'dberg and Ostrovskii concerning linear differential
equations with coefficients of completely regular growth,
\emph{Proc. Amer. Math. Soc.} \textbf{151} (2023), no.~5, 2097--2101.
\url{https://doi.org/10.1090/proc/16317}.

\bibitem{B1954}
R.~P.~Boas, Jr.,
\emph{Entire Functions},
Pure and Applied Mathematics, vol.~5,
Academic Press, New York, 1954.

\bibitem{c2023}
Y.-M.~Chiang, A.~Ching, and X.~Lin,
Holonomic Bessel modules and generating functions,
arXiv:2303.15496v1, 2023.
\url{https://arxiv.org/abs/2303.15496}.

\bibitem{C2009}
Y.-M.~Chiang and S.-J.~Feng,
On the growth of logarithmic differences, difference quotients
and logarithmic derivatives of meromorphic functions,
\emph{Trans. Amer. Math. Soc.} \textbf{361} (2009), no.~7, 3767--3791.
\url{https://doi.org/10.1090/S0002-9947-09-04663-7}.

\bibitem{C2016}
Y.-M.~Chiang and S.-J.~Feng,
On the growth of logarithmic difference of meromorphic functions
and a Wiman--Valiron estimate,
\emph{Constr. Approx.} \textbf{44} (2016), no.~3, 313--326.
\url{https://doi.org/10.1007/s00365-016-9324-8}.

\bibitem{Eremenko1998}
A.~Eremenko,
Extremal holomorphic curves for defect relations,
\emph{J. Anal. Math.} \textbf{74} (1998), no.~1, 307--323.
\url{https://doi.org/10.1007/BF02819454}.

\bibitem{G2026}
G.~G.~Gundersen and J.~Heittokangas,
\emph{Entire Functions: Classical and Modern Results, Special Functions
and Complex Differential Equations},
De Gruyter Textbook, De Gruyter, Berlin, 2026.
\url{https://doi.org/10.1515/9783112249130}.

\bibitem{H1994}
V.~P.~Havin and N.~K.~Nikolski (eds.),
\emph{Linear and Complex Analysis. Problem Book 3. Part II},
Lecture Notes in Mathematics, vol.~1574,
Springer-Verlag, Berlin, 1994.

\bibitem{H2015}
Janne~Heittokangas, Ilpo~Laine, Kazuya~Tohge and Zhi-Tao~Wen,
Completely regular growth solutions of second order complex
linear differential equations,
\emph{Ann. Acad. Sci. Fenn. Math.} \textbf{40} (2015), no.~2, 985--1003.
\url{https://doi.org/10.5186/aasfm.2015.4057}.

\bibitem{Hormander1990}
L.~H\"ormander,
\emph{The Analysis of Linear Partial Differential Operators~I:
Distribution Theory and Fourier Analysis},
2nd ed., Springer-Verlag, Berlin, 1990.

\bibitem{I2021}
K.~Ishizaki and Z.-T.~Wen,
Binomial series and complex difference equations,
\emph{J. Math. Anal. Appl.} \textbf{497} (2021), no.~1,
Paper No.~124844.
\url{https://doi.org/10.1016/j.jmaa.2020.124844}.

\bibitem{I2025}
K.~Ishizaki and Z.-T.~Wen,
Possible orders of transcendental meromorphic solutions
of linear difference equations with polynomial coefficients,
\emph{J. Anal. Math.} \textbf{155} (2025), no.~2, 459--483.
\url{https://doi.org/10.1007/s11854-025-0357-4}.

\bibitem{I2004}
K.~Ishizaki and N.~Yanagihara,
Wiman--Valiron method for difference equations,
\emph{Nagoya Math. J.} \textbf{175} (2004), 75--102.
\url{https://doi.org/10.1017/S0027763000008916}.

\bibitem{J1980}
B.~Ya.~Levin,
\emph{Distribution of Zeros of Entire Functions},
revised ed., Translations of Mathematical Monographs, vol.~5,
American Mathematical Society, Providence, RI, 1980.
Translated from the Russian by R.~P.~Boas, J.~M.~Danskin,
F.~M.~Goodspeed, J.~Korevaar, A.~L.~Shields, and H.~P.~Thielman.
\url{https://doi.org/10.1090/mmono/005}.

\bibitem{Levin1996}
B.~Ya.~Levin,
\emph{Lectures on Entire Functions},
Translations of Mathematical Monographs, vol.~150,
American Mathematical Society, Providence, RI, 1996.
\url{https://doi.org/10.1090/mmono/150}.

\bibitem{M1965}
A.~I.~Markushevich,
\emph{Theory of Functions of a Complex Variable}, vol.~II,
revised English ed., Prentice-Hall, Englewood Cliffs, NJ, 1965.
Translated and edited by R.~A.~Silverman.

\bibitem{P1984}
V.~P.~Petrenko,
\emph{Entire Curves},
Vishcha Shkola, Kharkov, 1984 (in Russian).

\bibitem{P1986}
C.~Praagman,
Fundamental solutions for meromorphic linear difference equations in the complex plane, and related problems,
\emph{J. Reine Angew. Math.} \textbf{369} (1986), 101--109.
\url{https://doi.org/10.1515/crll.1986.369.101}.

\bibitem{R1992}
L.~I.~Ronkin,
\emph{Functions of Completely Regular Growth},
Mathematics and its Applications (Soviet Series), vol.~81,
Kluwer Academic Publishers, Dordrecht, 1992.
Translated from the Russian by A.~Ronkin and I.~Yedvabnik.
\url{https://doi.org/10.1007/978-94-011-2418-8}.

\bibitem{S1989}
N.~Steinmetz,
Exceptional values of solutions of linear differential equations,
\emph{Math. Z.} \textbf{201} (1989), no.~3, 317--326.
\url{https://doi.org/10.1007/BF01214898}.

\bibitem{V1949}
G.~Valiron,
\emph{Lectures on the General Theory of Integral Functions},
Chelsea Publishing Company, New York, 1949.
Translated by E.~F.~Collingwood.

\bibitem{Wen2018}
Zhi-Tao~Wen, Gary~G.~Gundersen and Janne~Heittokangas,
Dual exponential polynomials and linear differential equations,
\emph{J. Differential Equations} \textbf{264} (2018), no.~1, 98--114.
\url{https://doi.org/10.1016/j.jde.2017.09.003}.
	\end{thebibliography}
\end{document}